\documentclass{gtart_h}

\def\ifplaintex{\expandafter\ifx\csname documentclass\endcsname\relax}

\def\gtp{{\mathsurround=0pt\it $\cal G\mskip-2mu$eometry \&\ 
$\cal T\!\!$opology $\cal P\!$ublications}}  

\def\recd{{\small Received:\qua\receiveddate\ifx\reviseddate\relax
\else\qquad Revised:\qua\reviseddate\fi\par}} 

\def\lognumber#1{\def\thelognumber{#1}}
\def\volumenumber#1{\def\thevolumenumber{#1}}
\def\volumeyear#1{\def\thevolumeyear{#1}}
\def\papernumber#1{\def\thepapernumber{#1}}
\def\pagenumbers#1#2{\def\startpage{#1}\def\finishpage{#2}}
\def\published#1{\def\publishdate{#1}}

\def\received#1{\def\receiveddate{#1}}

\def\accepted#1{\def\accepteddate{#1}}
\def\asciititle#1{\def\theasciititle{#1}}

\def\asciiaddress#1{\def\theasciiaddress{#1}}
\def\asciiemail#1{\def\theasciiemail{#1}}
\def\asciiurl#1{\def\theasciiurl{#1}}

\long\def\asciiabstract#1{\long\def\theasciiabstract{#1}}
\def\asciikeywords#1{\def\theasciikeywords{#1}}

\let\\\par\let\thelognumber\relax\let\thevolumenumber\relax
\let\thepapernumber\relax\let\thevolumeyear\relax\let\startpage\relax
\let\finishpage\relax\let\publishdate\relax\let\receiveddate\relax
\let\reviseddate\relax\let\accepteddate\relax\let\theasciititle\relax
\let\theasciiauthors\relax\let\theasciiaddress\relax
\let\theasciiabstract\relax\let\theasciikeywords\relax

\let\theasciiemail\relax
\let\theasciiurl\relax

\ifplaintex
\font\logobig=cmssbx10 scaled 3836
\font\logomed=cmssbx10 scaled 2557
\else
\font\logobig=cmssbx10 scaled 4200
\font\logomed=cmssbx10 scaled 2800
\fi

\long\def\makeagttitle{   
\count0=\startpage
\agt\hfill      
\hbox to 45truept{\vbox to 0pt{\vglue -13truept{\logomed A\kern -.37em{\logobig 
T}\kern -.38em G}\vss}\hss}
\break
{\small Volume \thevolumenumber\ (\thevolumeyear)
\startpage--\finishpage\nl
Published: \publishdate}

\vglue .25truein

{\parskip=0pt\leftskip 0pt plus
1fil\def\\{\par\smallskip}{\Large\bf\thetitle}\par\medskip} \vglue
0.05truein

{\parskip=0pt\leftskip 0pt plus 1fil\def\\{\par}{\sc\theauthors}
\par\medskip}%
 
\vglue 0.03truein 

{\small\leftskip 25truept\rightskip 25truept{\bf Abstract}\stdspace\theabstract

{\bf AMS Classification}\stdspace\theprimaryclass
\ifx\thesecondaryclass\relax\else; \thesecondaryclass\fi\par
{\bf Keywords}\stdspace \thekeywords\par}\vglue 7truept

}   

\ifplaintex
\font\phead=cmsl9 scaled 950
\font\pnum=cmbx10 scaled 913
\font\pfoot=cmsl9 scaled 950
\headline{\vbox to 0pt{\vskip -4.5mm\line{\small\phead\ifnum
\count0=\startpage ISSN 1472-2739 (on-line) 1472-2747 (printed)
\hfill {\pnum\folio}\else\ifodd\count0\def\\{ }%
\ifx\theshorttitle\relax\thetitle\else\theshorttitle\fi\hfill{\pnum\folio}
\else\def\\{ and }{\pnum\folio}\hfill\ifx\theshortauthors\relax\theauthors
\else\theshortauthors\fi\fi\fi}\vss}}
\footline{\vbox to 0pt{\vglue 0mm\line{\small\pfoot\ifnum\count0=\startpage
\copyright\ \gtp\hfill\else
\agt, Volume \thevolumenumber\ (\thevolumeyear)\hfill\fi}\vss}}
\else
\font=cmsl9 scaled 1050
\font\lnum=cmbx10 
\font=cmsl9 scaled 1050
\makeatletter
\def\@oddhead{{\small\lhead\ifnum\count0=\startpage ISSN 1472-2739 
(on-line) 1472-2747 (printed)\hfill {\lnum\number\count0}\else\ifodd\count0
\def\\{ }\ifx\theshorttitle\relax \thetitle \else\theshorttitle\fi\hfill
{\lnum\number\count0}\else\def\\{ and }{\lnum\number\count0}
\hfill\ifx\theshortauthors\relax 
\theauthors\else\theshortauthors\fi\fi\fi}}\def\@evenhead{\@oddhead}
\def\@oddfoot{\small\lfoot\ifnum\count0=\startpage\copyright\ \gtp\hfill\else
\agt, Volume \thevolumenumber\ (\thevolumeyear)\hfill\fi}
\def\@evenfoot{\@oddfoot}
\makeatother
\fi
\let\maketitlepage\makeagttitle
\let\makeshorttitle\maketitlepage
\let\maketitle\maketitlepage

\newwrite\gtoutfile
\long\gdef\makeheadfile{  
{\def\\{, }\def\s{ }
\immediate\openout\gtoutfile head.xxx
\immediate\write\gtoutfile{Proxy-for: \ifx\theasciiauthors\relax
\theauthors\else\theasciiauthors\fi\s<\ifx\theasciiemail\relax\theemail\else\theasciiemail\fi>}
\immediate\write\gtoutfile{\noexpand\\}
\immediate\write\gtoutfile{Authors: \ifx\theasciiauthors\relax
\theauthors\else\theasciiauthors\fi}
{\def\\{ }\immediate\write\gtoutfile{Title: \ifx\theasciititle\relax
\thetitle\else\theasciititle\fi}}
\immediate\write\gtoutfile{Subj-class: GT or SG, GR etc}
\immediate\write\gtoutfile{MSC-class: \theprimaryclass\ifx\thesecondaryclass\relax\else, \thesecondaryclass\fi}
\immediate\write\gtoutfile{Journal-ref: Algebr. Geom. Topol. \thevolumenumber\s
(\thevolumeyear) \startpage-\finishpage}
\immediate\write\gtoutfile{Comments: Published by Algebraic and
Geometric Topology at}
\immediate\write\gtoutfile{\s\s\s  http://www.maths.warwick.ac.uk/agt/AGTVol\thevolumenumber/agt-\thevolumenumber-\thepapernumber.abs.html}
\immediate\write\gtoutfile{\noexpand\\}
\immediate\write\gtoutfile{}
\ifx\theasciiabstract\relax
\immediate\write\gtoutfile{\theabstract}\else
\immediate\write\gtoutfile{\theasciiabstract}\fi
\immediate\write\gtoutfile{}
\immediate\write\gtoutfile{\noexpand\\}
\immediate\write\gtoutfile{}
\immediate\closeout\gtoutfile}}  

\def\maketitlepage{\makeagttitle\makeheadfile}
\let\makeshorttitle\maketitlepage
\let\maketitle\maketitlepage

\lognumber{50}
\volumenumber{4}
\volumeyear{2004}
\papernumber{50}
\pagenumbers{1145}{1153}
\received{13 June 2004} 
\accepted{16 September 2004}
\published{1 December 2004}

\usepackage{amsmath,amssymb}

\newcommand{\R}{{\mathbb R}}
\newcommand{\C}{{\mathbb C}}
\newcommand{\Z}{{\mathbb Z}}
\newcommand{\Q}{{\mathbb Q}}
\renewcommand{\H}{{\mathbb H}}
\newcommand{\maps}{\colon\thinspace}
\DeclareMathOperator{\tr}{tr}
\newcommand{\PSL}[2]{\mathrm{PSL}_{#1} #2}
\newcommand{\SU}[1]{\mathrm{SU}_{ #1}}
\newcommand{\SL}[2]{\mathrm{SL}_{#1} #2}
\newcommand{\abs}[1]{{\left| #1 \right|}}

\theoremstyle{plain} 

\swapnumbers
\newtheorem{theorem}{Theorem}[section]

\newtheorem{lemma}[theorem]{Lemma}
\newtheorem{corollary}[theorem]{Corollary}

\newtheorem{claim}[theorem]{Claim}

\theoremstyle{definition}

\theoremstyle{remark}
\newtheorem{rmk}[theorem]{Remark}

\makeatletter
  \let\c@theorem=\c@subsection
  \let\p@theorem=\p@subsection
  
  \let\cl@theorem=\cl@subsection
   \let\c@figure=\c@subsection
  \let\p@figure=\p@subsection
  
  \let\cl@figure=\cl@subsection
  \let\c@table=\c@subsection
  \let\p@table=\p@subsection
  
  \let\cl@table=\cl@subsection
  \let\c@equation=\c@subsection
  \let\p@equation=\p@subsection
  
  \let\cl@equation=\cl@subsection
\makeatother

\begin{document}

\title{Non-triviality of the $A$-polynomial for knots in $S^3$}
\asciititle{Non-triviality of the A-polynomial for knots in S^3}

\authors{Nathan M. Dunfield\\Stavros Garoufalidis}
\address{Mathematics 253-37, California Institute of 
Technology\\Pasadena,  CA 91125,  USA}

\secondaddress{School of Mathematics, Georgia Institute of 
Technology\\Atlanta, GA 30332-0160, USA }

\asciiaddress{Mathematics 253-37, California Institute of 
Technology\\Pasadena,  CA 91125,  USA\\and\\School of Mathematics, 
Georgia Institute of Technology\\Atlanta, GA 30332-0160, USA }

\gtemail{\mailto{dunfield@caltech.edu}\qua{\rm 
and}\qua\mailto{stavros@math.gatech.edu}}
\gturl{\url{http://www.its.caltech.edu/~dunfield}\qua{\rm 
and}\qua{http://www.math.gatech.edu/~stavros}}

\asciiemail{dunfield@caltech.edu, stavros@math.gatech.edu}
\asciiurl{http://www.its.caltech.edu/ dunfield,
http://www.math.gatech.edu/ stavros }

\begin{abstract} 
  The $A$-polynomial of a knot in $S^3$ defines a complex plane curve
  associated to the set of representations of the fundamental group of
  the knot exterior into $\SL{2}{\C}$.  Here, we show that a
  non-trivial knot in $S^3$ has a non-trivial $A$-polynomial.  We
  deduce this from the gauge-theoretic work of Kronheimer and Mrowka
  on $\SU{2}$-representations of Dehn surgeries on knots in $S^3$.  As
  a corollary, we show that if a conjecture connecting the colored
  Jones polynomials to the $A$-polynomial holds, then the colored
  Jones polynomials distinguish the unknot.
\end{abstract}

\asciiabstract{%
  The A-polynomial of a knot in S^3 defines a complex plane curve
  associated to the set of representations of the fundamental group of
  the knot exterior into SL(2,C).  Here, we show that a
  non-trivial knot in S^3 has a non-trivial A-polynomial.  We
  deduce this from the gauge-theoretic work of Kronheimer and Mrowka
  on SU_2-representations of Dehn surgeries on knots in S^3.  As
  a corollary, we show that if a conjecture connecting the colored
  Jones polynomials to the A-polynomial holds, then the colored
  Jones polynomials distinguish the unknot}

\asciikeywords{Knot, A-polynomial, character variety,  Jones polynomial}
\maketitle

\section{Introduction}
\label{sec:intro}

Roughly speaking, the $A$-polynomial of a knot $K$ in $S^3$ describes
the $\SL 2 \C$-representations of the knot complement, as viewed from
the boundary.  In a little more detail, let $M$ be the exterior of $K$.
The boundary of $M$ is a torus, whose fundamental group $\pi_1(\partial M) =
\Z^2$ comes with a natural meridian-longitude basis $(\mu, \lambda)$.
Consider a representation $\rho \maps \pi_1(M) \to \SL 2 \C$  The
restriction of $\rho$ to $\pi_1(\partial M)$ has a simple form, since a pair
of commuting 2-by-2 matrices are typically simultaneously
diagonalizable, i.e.~$\rho$ can be conjugated so that:
\[
\rho(\mu)=\begin{pmatrix}
M & 0 \\ 0 & M^{-1}
\end{pmatrix} 
\qquad \mbox{and} \qquad 
\rho(\lambda )=\begin{pmatrix}
L & 0 \\ 0 & L^{-1}
\end{pmatrix}. 
\]
The possible eigenvalues $(M, L)$ of $(\rho(\mu), \rho(\lambda))$ as $\rho$
varies form an complex algebraic subvariety of $\C^2$.  The
$A$-polynomial is the defining equation for the $1$-dimensional part
of this subvariety; that is, it describes a plane curve whose points
correspond to the restrictions of representations to
$\pi_1(\partial M)$.

The $A$-polynomial of a knot, which was introduced by Cooper et al.~in
\cite{CCGLS}, has deep connections to the topology and geometry of
$M$  As the group of isometries of hyperbolic 3-space is
$\PSL{2}{\C}$, the $A$-polynomial is connected to the study of
deformations of (incomplete) hyperbolic structures on $M$.  For
example, the variation of the volume of hyperbolic structures on $M$
depends only on their restriction to the boundary torus, and is controlled
entirely by the $A$-polynomial.  On the topological side, the sides of
the Newton polygon of the $A$-polynomial give rise to incompressible
surfaces in $M$.

Here, we address the basic question: can $A$-polynomial distinguish
the unknot from all other knots in $S^3$?  The $A$-polynomial of the
unknot is simply $L - 1$.  The $A$-pol\-y\-no\-mi\-al always contains a
factor of $L - 1$ coming from reducible representations; we say that
the $A$-polynomial is non-trivial if it has an additional factor.
Perhaps for some non-trivial knots, there are no other
representations, or they don't deform in ways that change the holonomy
on the boundary.  Our main result shows that this does not happen, and
hence the $A$-polynomial distinguishes the unknot: 
\begin{theorem}\label{main-thm}
  A non-trivial knot in $S^3$ has a non-trivial $A$-polynomial.
  Moreover, the $A$-pol\-y\-no\-mi\-al is not a power of $L - 1$.
\end{theorem}
Steve Boyer and Xingru Zhang independently proved
Theorem~\ref{main-thm} using a similar approach \cite{BoyerZhang}.

We deduce Theorem~\ref{main-thm} as a direct corollary of the
following deep theorem of Kronheimer and Mrowka:
\begin{theorem}{\rm\cite{KM}}\label{KM-thm}\qua
  Let $K$ be a non-trivial knot in $S^3$.  For $r \in \Q$, let $M_r$
  be the 3-manifold which is the $r$ Dehn surgery on $K$.  If $\abs{r}
  \leq 2$, then there exists a homomorphism $\pi_1(M_r) \to \SU{2}$ with
  non-cyclic image.
\end{theorem}
Their proof uses gauge theory; in addition to their own major
contributions, the proof relies on Gabai's theorem that the
zero-surgery on knot has a taut-foliation, Eliashberg and Thurston's
work connecting foliations to contact structures, Eliashberg's proof
that contact 3-manifolds embed in symplectic
4\nobreakdash-\hspace{0pt}manifolds, Taubes' non-vanishing theorem for
Seiberg-Witten invariants of symplectic 4-manifolds, and Feehan and
Leness' work connecting the Seiberg-Witten and Donaldson invariants.

Theorem~\ref{main-thm} was previously known for all non-satellite
knots for simple geometric reasons, as we now describe.  When $M$ is
hyperbolic, we have the holonomy representation $\pi_1(M) \to \SL 2 \C$
of the complete hyperbolic structure; Thurston showed in his
Hyperbolic Dehn Surgery Theorem that this representation has a complex
curve of deformations which change the holonomy along the boundary
\cite{Th}.  Thus, in this case, the $A$-polynomial is non-trivial.
Non-hyperbolic knots are torus knots or satellites. For torus knots, a
simple calculation shows they have non-trivial $A$-polynomial
\cite{CCGLS}.  Satellite knots are those which have closed
incompressible tori in their complements.  One can look at the
resulting geometric decomposition, and try to understand how the
representations of each piece could glue together to give a
representation of all of $\pi_1(M)$; however, this seems quite
difficult to do in general.

Since our proof of Theorem~\ref{main-thm} is based on the existence of
$\SU{2}$ representations, we really show that if one looks only at
representations $\rho \maps \pi_1(M) \to \SU 2$, then the eigenvalues
$(M, L)$ of $(\rho(\mu), \rho(\lambda))$ sweep out a real 1-dimensional
subset of the unit torus in $\C^* \times \C^*$.  This is interesting even in the case
of hyperbolic knots.

\subsection{Connection to the Jones polynomial}

While the $A$-polynomial arose from the study of hyperbolic geometry,
it turns out to have connections to seemingly disparate parts of
low-dimensional topology, including the Jones
polynomial.  As we will now explain, the non-triviality of the
$A$-polynomial of a knot has implications to the strength of the
colored Jones function. The latter is essentially the sequence of
Jones polynomials of a knot and its connected parallels.  In
\cite{GL}, it was proven that the colored Jones function of a knot is
a sequence of Laurent polynomials which satisfy a $q$-difference
equation.  It was observed by the second author in \cite{Ga} that one
can choose the $q$-difference equation in a canonical manner.  The
corresponding operator to this $q$-difference equation is an element
of the non-commutative ring
$$
\Z[q^{\pm}]\langle Q^{\pm}, E^{\pm} \rangle /( EQ-qQE)
$$
of Laurent polynomials in $E$ and $Q$ that satisfy the commutation
relation $EQ=qQE$. 

This operator defines the so-called non-commutative $A$-polynomial of
a knot.  In \cite{Ga}, the second author conjectured that specializing
the non-commutative $A$-polynomial at $q=1$ coincides with the
$A$-polynomial of a knot after the change of variables $(E,Q)=(L,M^2)$
(there may also be changes in the multiplicities of factors and
polynomials in $Q$).  This is called the AJ Conjecture, and an
immediate consequence of Theorem~\ref{main-thm} is:
\begin{corollary}
If the AJ Conjecture holds, then the colored Jones function distinguishes the
unknot.
\end{corollary}

\subsection{Connection to contact homology}  

Another surprise is that the $A$-polynomial is connected with contact
geometry.  Consider the unit conormal bundle to $\R^3$, denoted
$ST^*(\R^3)$, which has a natural contact structure.  If $K$ is a knot
in $\R^3$ then the unit conormal bundle to $K$ is a Legendrian 2-torus
$L$ inside $ST^*(\R^3)$.  Lenny Ng has constructed a homology theory
for knots in $S^3$, the \emph{framed knot contact homology}, which is
strongly believed to be Eliashberg-Hofer of contact homology of the
pair $(ST^*(\R^3), L)$ \cite{Ng}.  Ng has shown that the
$A$-polynomial can be derived from the simplest piece of the framed
knot contact homology.  Combining this with Theorem~\ref{main-thm}, he
proves:

\begin{theorem}{\rm\cite[Prop.~5.9]{Ng}}\qua
  The framed knot contact homology distinguishes the unknot from any other
  knot in $S^3$.
\end{theorem}

\subsection{Acknowledgments}  

Both authors were partially supported by the U.S. National Science
Foundation, and Dunfield was also partially supported by the Sloan Foundation.  

\section{Proofs}

We begin by reviewing the definition of the $A$-polynomial for a
compact 3-manifold $M$ whose boundary is a torus (for details, see
\cite{Sh, CCGLS}).  Let $R(M)$ denote the set of representations
$\pi_1(M) \to \SL 2 \C$, which is an affine algebraic variety over $\C$.
It is natural to study such representations up to inner automorphisms
of $\SL 2 \C$, so consider the \emph{character variety}, $X(M)$, which
is the quotient of $R(M)$ under the action of $\SL 2 \C$ by
conjugation.  Technically, one has to take the algebro-geometric
quotient to deal with orbits of reducible representations which are
not closed; in this way $X(M)$ is also an affine complex algebraic
variety.

To define the $A$-polynomial, we first need to understand the character
variety $X(\partial M)$ of the torus $\partial M$.  The fundamental group of $\partial
M$ is just $\Z \times \Z$, and fix generators $(\mu, \lambda)$.  Since $\pi_1(\partial
M)$ is commutative, any representation $\rho \maps \pi_1(\partial M) \to \SL 2
\C$ is reducible, that is, has a global fixed point for the M\"obius
action on $P^1(\C)$.  Moreover, if no element of $\rho(\pi_1(\partial M))$ is
parabolic, $\rho$ is conjugate to a diagonal representation with
\[
\rho(\mu)=\begin{pmatrix}
M & 0 \\ 0 & M^{-1}
\end{pmatrix} 
\qquad \mbox{and} \qquad 
\rho(\lambda )=\begin{pmatrix}
L & 0 \\ 0 & L^{-1}
\end{pmatrix}. 
\]
As such, $X(\partial M)$ is approximately $\C^* \times \C^*$ with
coordinates being the eigenvalues $(M,L)$.  This isn't quite right, as
switching $(M, L)$ with $(M^{-1}, L^{-1})$ gives a conjugate
representation.  In fact, $X( \partial M)$ is exactly the quotient of
$\C^* \times \C^*$ under the involution $(M, L) \mapsto (M^{-1},
L^{-1})$.

Now the inclusion $i \maps \partial M \to M$ induces a regular map $i^* \maps
X(M) \to X(\partial M)$ via restriction of representations from $\pi_1(M)$ to
$\pi_1(\partial M)$.  Let $V$ be the (complex) 1-dimensional part of
$i^*\big(X(M)\big)$.  More precisely, take $V$ to be the union of the
1-dimensional $i^*(X)$, where $X$ is an irreducible component of
$X(M)$.  The curve $V$ is used to define the $A$-polynomial.  To
simplify things, we look at the plane curve $\overline{V}(M)$ which is
inverse image of $V$ under the quotient map $\C^* \times \C^* \to X(\partial M)$
The $A$-polynomial is the defining equation for $\overline{V}(M)$; it
is a polynomial in the variables $M, L$. Since all the maps involved
are defined over $\Q$, the $A$-polynomial can be normalized to have
integral coefficients.

In the definition of the $A$-polynomial, we looked only at those
irreducible components where $i^*(X)$ is 1-dimensional.  In the proof
of Theorem~\ref{main-thm}, we will need the following standard lemma
to show that we do not overshoot our goal of showing that the
$A$-polynomial is non-trivial.
\begin{lemma}\label{lemma}
  Let $X$ be an irreducible component of $X(M)$.  Then $i^*(M)$ has
  dimension $0$ or $1$.
\end{lemma}
\begin{proof}
  There are two proofs of this in the literature, and we include
  sketches of both to make this paper more self-contained.
  
  To prove the lemma, we just need to rule out the possibility that
  $i^*(X)$ is 2-dimensional, and thus a Zariski-open subset of $X(\partial
  M)$.  The approach in \cite[\S4.5]{CCGLS} is to introduce the notion
  of the volume of a representation $\rho \maps \pi_1(M) \to \SL{2}{\C}$
  (see also \cite[\S2.5]{Dunfield} and \cite[\S4]{Franc} for a more
  complete definition of the volume).  This gives a natural function
  $\mathrm{Vol} \maps X(M) \to \R$.  Then Schl\"afli's formula for the change
  of volume of a family of polyhedra in $\H^3$ shows that the
  derivative of $\mathrm{Vol}$ depends only on the restriction of
  representations to $\pi_1(\partial M)$.  This leads to a $1$-form on $X(\partial
  M)$ which must be exact on $i^*(X(M))$.  This form is not exact on
  any Zariski-open subset of $X(\partial M)$, and hence $i^*(X)$ is at most
  $1$-dimensional.
  
  The other argument is to observe that if $i^*(X)$ were
  2-dimensional, it would let us construct ideal points of $X(M)$
  where the associated surface has whatever boundary slope we want.
  This would contradict Hatcher's finiteness theorem for boundary
  slopes.  In more detail, start with a slope $\alpha \in \pi_1(\partial M)$ and
  let $\beta$ be a complementary slope.  Choose a $c \in \C$ so that the
  curve $Y$ in $X(\partial M)$ given by $\tr_\alpha = c$ has $i^*(X) \cap Y$ dense
  in $Y$.  Choose a curve $\widetilde{Y} \subset X$ whose image under $i^*$
  is dense in $Y$.  As $\tr_\beta$ is non-constant on $Y$, there is an
  ideal point $p$ of $\widetilde{Y}$ where $\tr_\beta$ has a pole.  Since
  $\tr_\alpha$ is constant on $Y$, an incompressible surface associated to
  the ideal point $p$ must have boundary slope $\alpha$.  But Hatcher
  showed that there are only finitely many $\alpha$ which are boundary
  slopes of incompressible surfaces \cite{Hatcher}, a contradiction.
\end{proof}
  
  When $M$ is the exterior of a knot in $S^3$, then, up to orientation
  conventions, there is a canonical meridian-longitude basis $(\mu,
  \lambda)$ for $\pi_1(\partial M)$, and one uses this basis when writing the
  $A$-polynomial.  Since we are interested in the non-triviality of
  the $A$-polynomial, we need to discuss the conventions for dealing
  with the reducible representations.  When $M$ is the exterior of a
  knot in $S^3$, one has $H_1(M, \Z) = \Z$, and so there are many
  reducible representations which factor via: $\pi_1(M) \to \Z \to \SL 2
  \C$.  Irreducible components of $X(M)$ either consist solely of
  reducible representations, or have a Zariski-open subset of
  irreducible representations.  In the case of the exterior of a knot in
  $S^3$, there is a single irreducible component of $X(M)$ consisting
  entirely of reducible representations.  This component contributes a
  factor of $L - 1$ to the $A$-polynomial. Some authors exclude this
  factor from the $A$-polynomial, and define the curve $V$ above to be
  the image under $i^*$ of the irreducible components of $X(M)$ which
  contain an irreducible representation.  To say the $A$-polynomial is
  non-trivial, we mean that it does not just consist of the $L - 1$
  coming from the reducible representations.  We will now show that
  the $A$-polynomial of a non-trivial knot in $S^3$ is non-trivial,
  and, moreover, is not just a power of $L - 1$.

\begin{proof}[Proof of Theorem~\ref{main-thm}]
  Let $M$ be the exterior of a non-trivial knot in $S^3$.  Let $X'(M)$
  denote $X(M)$ minus the component consisting of reducible
  representations, and let $V'$ be the union of the 1-dimensional
  $i^*(X)$ where $X$ is an irreducible component of $X'(M)$.  The main
  part of the theorem is that $V'$ is non-empty.  To this end, we will
  show:

  \begin{claim}  \label{main-claim}
    There exists  an infinite collection of
    irreducible representations 
    $\rho_n \maps \pi_1(M) \to \SL 2 \C$
    whose restrictions to $\pi_1(\partial M)$ are all distinct in $X(\partial M)$.  
  \end{claim}
  
  Before proving the claim, let us deduce $V' \neq \emptyset$ from it.
  Assuming the claim, then as a $0$-dimensional algebraic variety
  consists of finitely many points, there must be some irreducible $X$
  in $X'(M)$ so that the dimension of $i^*(X)$ is at least 1.  By
  Lemma~\ref{lemma}, the dimension of $i^*(X)$ must be exactly one, and
  so $V' \neq \emptyset$.
  
  To prove the claim, we use the $\SU 2$ representations given by
  Theorem~\ref{KM-thm}.  Let $M_{1/n}$ denote the $1/n$-filling of
  $M$.  By Theorem~\ref{KM-thm}, for each non-zero $n \in \Z$ we have
  a representation $\rho_n \maps \pi_1(M_{1/n}) \to \SU 2$ with non-cyclic
  image.  First, we claim that the $\rho_n$ are irreducible as
  representations into the larger group $\SL 2 \C$.  Suppose $\rho_n$
  were reducible.  Since $H_1(M_{1/n}, \Z) = 0$, the group $G =
  \pi_1(M_{1/n})$ satisfies $G = [G, G]$.  As $\rho_n$ is reducible, and
  commutators of elements of $\SL 2 \C$ with a common fixed point are
  parabolic with trace $2$, it follows that $\tr( \rho_n(\gamma) ) = 2$ for
  all $\gamma \in G$.  But the only element of $\SU 2$ with trace 2 is the
  identity, and so $\rho_n$ would be trivial, a contradiction.  So
  $\rho_n$ is irreducible.

  As $\pi_1(M_{1/n})$ is a quotient of $\pi_1(M)$, we will regard $\rho_n$
  as a representation of $\pi_1(M)$ into $\SU 2 \leq \SL 2 \C$.  To prove
  Claim~\ref{main-claim}, we need to show that the restrictions of the
  $\rho_n$ to $\pi_1(\partial M)$ gives an infinite collection of points in
  $X(\partial M)$.  Two representations of $\pi_1(\partial M)$ into $\SU 2$ which
  correspond to the same point in $X(\partial M)$ are actually
  conjugate---this is because they both must be conjugate to diagonal
  representations (this isn't quite true for $\SL 2 \C$, where
  distinct parabolic representations get amalgamated).  Because of
  this, to prove the Claim~\ref{main-claim} it suffices to show that
  the kernels $K_n$ of the $\rho_n$ give an infinite collection of
  distinct subgroups of $\pi_1(\partial M) = \Z^2$.
 
  For $\alpha$ a slope in $\partial M$, note that $\rho_n$ extends to $\pi_1(M_\alpha)$
  if and only if $\alpha \in K_n$.  As $\rho_n$ comes from $M_{1/n}$, we have
  $(1,n) \in K_n$ for each $n \neq  0$.  As the $1/0$ filling gives $S^3$,
  which is simply connected, we have $(1,0) \notin K_n$.  Because of this,
  Claim~\ref{main-claim} follows from directly from the following
  lemma with $\gamma$ the line $x = 1$:

  \begin{lemma}
    \label{lem.2}
    Suppose $\gamma$ is a line in $\R^2$ which contains infinitely many
    lattice points of $\Z^2$, and which does not contain $0$.
    Consider a collection $K_n$ of subgroups of $\Z^2$ whose union,
    $K$, contains all but finitely many of the lattice points on $\gamma$.
    Suppose, in addition, that there is a lattice point on $\gamma$
    which is not in $K$.  Then there are infinitely many distinct
    $K_n$.
  \end{lemma}

\begin{proof}
  Assume that there are finitely many $K_n$.  If $K_n$ has rank less
  than 2, then $K_n$ is contained in a line through the origin, and so
  $K_n$ intersects $\gamma$ in at most one point.  So we can throw out
  all of the $K_n$ of rank less than 2, and still have $\gamma - K$ finite.
  
  So we can assume that $\Z^2 / K_n$ is finite for each $n$.  Let $L$
  be the intersection of $K_n$; as there are finitely many $K_n$, the
  subgroup $L$ is also a finite-index subgroup of $\Z^2$.  Now let
  $\gamma'$ be the line parallel to $\gamma$ which passes through the
  origin.  As $\Z^2/L$ is finite, the subgroup $H = \gamma' \cap L$ is
  infinite.  Let $v_0$ be the given point in $\gamma \setminus K$.  Then if $h \in
  H$, we have that $v_0 + h$ is also in $\gamma \setminus K$ since if $v_0 + h$
  is in some $K_n$, then so is $v_0 = (v_0 + h) - h$.  But $H$ is
  infinite, and thus so is $\{ v_0 + h \}$, which contradicts that $\gamma
  \setminus K$ is finite.  Thus we must have an infinite collection of
  distinct $K_n$.
\end{proof}

To complete the proof of Theorem~\ref{main-thm}, we need to show that
the $A$-polynomial of $M$ is not a power of $L - 1$.  Assume the
contrary  Consider the point $p_n = (m_n, l_n) \in \C^* \times \C^*$
corresponding to the restriction of the representation $\rho_n$ to
$\pi_1(\partial M)$.  As $\rho_n$ comes from the $(1, n)$ filling of $M$, we
have that $m_n l_n^n = 1$.  By the above argument, all but finitely
many of the pairs $(m_n, l_n)$ satisfy the $A$-polynomial, and hence
$l_n = 1$.  Then for such $n$, the relation $m_n l_n^n = 1$ implies
that $m_n = 1$.  As $\rho_n$ has image in $\SU{2}$, this implies that
$\rho_n$ is trivial when restricted to $\pi_1(\partial M)$. But then $\rho_n$
factors over to the $S^3$ surgery, a contradiction.  Thus the
$A$-polynomial is not a power of $L - 1$.
\end{proof}

\begin{rmk}
  Lemma~\ref{lem.2} has other applications to studying Dehn filling.
  For instance, consider a non-trivial knot $K$ in $S^3$ with exterior
  $M$.  In relation to the Virtual Haken Conjecture, this lemma
  implies there is a infinite sequence $n_k$ of non-zero integers so
  that the degree of the smallest non-trivial cover of $M_{1/n_k}$
  goes to infinity as $k \to \infty$.
\end{rmk}

\Addresses\recd

\end{document}